\documentclass[11pt]{amsart}

\usepackage{latexsym}
\usepackage{amsmath}
\usepackage{amssymb}
\usepackage{amscd}
\usepackage{multicol}
\usepackage{amsfonts}
\usepackage{amsthm}
\usepackage{psfrag}
\usepackage{epsfig}

\newtheorem{thm}{Theorem}
\newtheorem{theorem}[thm]{Theorem}

\newtheorem{lemma}[thm]{Lemma}
\newtheorem{prop}[thm]{Proposition}

\newtheorem{example}[thm]{Example}
\newtheorem*{exm*}{Example}

\numberwithin{thm}{section}

\newcommand{\ft}{{\mathfrak{t}}}

\newcommand{\fg}{{\mathfrak{g}}}

\renewcommand{\SS}{{\mathbb{S}}}
\newcommand{\RR}{{\mathbb{R}}}
\newcommand{\ZZ}{{\mathbb{Z}}}
\newcommand{\QQ}{{\mathbb{Q}}}
\newcommand{\CC}{{\mathbb{C}}}

\newcommand{\cA}{{\mathcal{A}}}

\newcommand{\cS}{{\mathcal{S}}}

\newcommand{\cC}{{\mathcal{C}}}
\newcommand{\cR}{{\mathcal{R}}}
\newcommand{\cO}{{\mathcal{O}}}

\DeclareMathOperator{\Inv}{Inv}

\title[Moment Map Degeneration]{Positivity of Equivariant Schubert Classes Through Moment Map Degeneration}

\author[Catalin Zara]{Catalin Zara}
\address{Department of Mathematics, UMass Boston, 100 Morrissey Blvd., Boston, MA 02125}

\begin{document}

\maketitle

\begin{abstract}
For a flag manifold $M=G/B$ with the canonical torus action, the $T-$equivariant cohomology is generated by
equivariant Schubert classes, with one class $\tau_u$ for every element $u$ of the Weyl group $W$.
These classes are determined by their restrictions to the fixed point set $M^T \simeq W$, and the restrictions
are polynomials with nonnegative integer coefficients in the simple roots.

The main result of this article is a positive formula for computing $\tau_u(v)$ in types A, B, and C. To obtain this formula we identify $G/B$ with a generic co-adjoint orbit and use a result of Goldin and Tolman to compute $\tau_u(v)$ in terms of the induced moment map. Our formula, given as a sum of contributions of certain maximal ascending chains from $u$ to $v$, follows from a systematic degeneration of the moment map, corresponding to degenerating the co-adjoint orbit. In type A we prove that our formula is manifestly equivalent to the formula announced by Billey in \cite{Bi}, but in type C, the two formulas are not equivalent.
\end{abstract}

\tableofcontents

\section{Introduction}

Let $G$ be a connected, complex, semisimple Lie group, $B \subset G$ a Borel subgroup, corresponding to a set of simple roots $\{ \alpha_1,\ldots,\alpha_n\}$. Let $T^{\CC} \subset B$ be a maximal torus, and $T \subset T^{\CC}$ a compact real form of $T^{\CC}$. The torus $T$ acts on the flag manifold $M=G/B$ by left multiplication on $G$, and the fixed point set $M^T$ corresponds bijectively with the Weyl group $W$.

The equivariant cohomology ring $H_T^*(M)= H_T^*(M,\QQ)$ is a free module over $H_T^*(pt,\QQ) \simeq \QQ[\alpha_1,\ldots,\alpha_{n}]$, and that is also true for integral, and not just rational coefficients. A basis of this module is given by equivariant Schubert classes, $\{ \tau_u \}_{u\in W}$, with a class $\tau_u$ for each fixed point $u \in W = M^T$. The equivariant Schubert class $\tau_u \in H_T^*(M)$ is the class induced by the $T-$equivariant cycle $X_u = \overline{B_{-}uB}/B$, where $B_{-}$ is the opposite Borel subgroup. The pull-back map $H_T^*(M) \to H_T^*(M^T)$ is injective, so an equivariant class is determined by its values on $M^T$.

In particular, $\tau_u$ is determined by the values $\tau_u(v) \in \QQ[\alpha_1,\ldots,\alpha_{n}]$, for $v \in W$. The integral positivity property of the equivariant Schubert class $\tau_u$ states that for all $v \in W$ we have
$$
  \tau_u(v) \in \mathbb{Z}_{\geqslant 0} [\alpha_1,\ldots,\alpha_{n}].
$$
In other words, $\tau_u(v)$ can be written as a
polynomial with {\em nonnegative integer} coefficients in the simple roots
$\alpha_1$,...,$\alpha_{n}$. Billey \cite{Bi} gives an \emph{integral positive formula} for $\tau_u(v)$:
a formula for $\tau_u(v)$ as a sum of polynomials in
$\alpha_1$, ..., $\alpha_{n}$ with non-negative integer coefficients.

The main result of this article, formulated in Theorem~\ref{thm:main-general}, is a new positive formula for computing $\tau_u(v)$ in types $A$, $B$, and $C$. To obtain this formula we identify $G/B$ with a generic co-adjoint orbit and use a recent result of Goldin and Tolman (\cite{GT}) to compute $\tau_u(v)$ in terms of the induced moment map. Our formula, given as a sum of contributions of certain maximal ascending chains from $u$ to $v$, follows from a systematic degeneration of the moment map, corresponding to degenerating the co-adjoint orbit. The resulting formula is integral in types $A$ and $C$, but only rational in type $B$ (and for $G_2$). In type $A$ our formula is an alternative formulation of the positive announced by Billey, but in type $C$ the two formulas, while giving the same answer, are not manifestly equivalent.

The main methods in this paper involve the combinatorics of Weyl groups. The type $A$ formula has been circulated earlier, and motivated a separate project by Sabatini and Tolman (\cite{ST}). Their geometrical methods allowed them to independently obtain, as applications, the same formulas for types $B$ and $C$, but from a different perspective.

The content of the paper is organized as follows. In Section~\ref{sec:combschub} we give a combinatorial description of $\tau_u$ and in Section~\ref{sec:subwordformula} we recall Billey's positive formula, that uses subwords of a reduced word for $v$. In Section~\ref{sec:GTformula} we apply the result of Goldin and Tolman to co-adjoint orbits and obtain a formula for $\tau_u(v)$ in terms of maximal ascending chains. This formula depends on the particular generic co-adjoint orbit and in Section~\ref{sec:limits} we use moment map/orbit degeneration to get a simpler formula, in which only some chains have non-zero contributions. In Section~\ref{sec:chainformula} we prove that the simpler formula is a positive formula, integral in types $A$ and $C$ and rational in type $B$, and in Section~\ref{sec:TypeAchains} we give an explicit version in type $A$. We prove the equivalence of the chain formula with the subword formula in Section~\ref{sec:equivalence}, after first constructing, in Section~\ref{sec:chains-subwords}, a general map connecting chains and subwords.

I would like to thank Rebecca Goldin, Sue Tolman, Silvia Sabatini, and Allen Knutson for explaining their work, and to all of the above and Victor Guillemin for many stimulating discussions and valuable comments while writing this paper. Many thanks are also due to John Stembridge for writing the \emph{Maple} packages \verb"coxeter" and \verb"weyl", which have been used to verify several conjectured formulas before proving them.

\section{Combinatorial Schubert Classes}
\label{sec:combschub}

The injective morphism $H_T^*(M) \to H_T^*(M^T)$ identifies $H_T^*(M)$ with a subring of $H_T^*(M^T) = \text{Maps}(W,\SS)$, where $\SS = \QQ[\alpha_1,\ldots,\alpha_{n}]$. Not all maps represent classes in $H_T^*(M)$, and the ones that do can be described using a discrete structure involving a regular graph $\Gamma=(V,E)$, and a labeling of the oriented edge of $\Gamma$ by elements of $\ft^*$, the dual of the Lie algebra of $T$.

The vertices $V$ of $\Gamma$ correspond to the Weyl group $W$. Two vertices $u$ and $v$ are joined by an edge if and only if they differ by a reflection $s_{\beta}$, for some positive root $\beta$. Note that $us_{\beta} = s_{u\beta}u$, so if $u$ and $v$ differ by a reflection to the left, they also differ by a reflection to the right.  If $v=us_{\beta}$, with $\beta \succ 0$ a positive root, then the oriented edge $e=(u,v)$ of $\Gamma$ is labeled by $\alpha(u,v) =u\beta$ and is called \emph{ascending} if $\alpha(u,v) \succ 0$. We will also use $\alpha(e)$, $\alpha_e$, or $\alpha_{u,v}$ for $\alpha(u,v)$. The pair $(\Gamma, \alpha)$ is called the \emph{GKM graph} of $(M,T)$.

An assignment $f \colon W \to \SS = \QQ[\alpha_1,\ldots,\alpha_{n}]$ is a \emph{cohomology class} if
\begin{equation*}
  f(v) - f(u) \in \alpha(e) \SS \, ,
\end{equation*}
for every edge $e=[u,v]$ of $\Gamma$. The class $f$ is called homogeneous of degree $k$ if, for every
$u \in W$, the polynomial $f(u)$ is homogeneous of degree $k$.
The \emph{cohomology ring} $H_{\alpha}^*(\Gamma)$ is the graded subring of
$\text{Maps}(W,\mathbb{S})$ consisting of all classes. Then
$$H_{T}^{2k}(M) \simeq H_{\alpha}^k(\Gamma) \; ,$$
and since $H_T^{odd}(M) = 0$, that means that $H_T^*(M) \simeq H_{\alpha}^*(\Gamma)$. We will in general identify a class $f \in H_T^*(M)$ with its image in $H_{\alpha}^*(\Gamma)$. In particular, we will denote by $\tau_u$ both the equivariant Schubert class in $H_T^*(M)$ and its image in $H_T^*(\Gamma)$. Sometimes we will refer to $\tau_u \in H_T^*(\Gamma)$ as the \emph{combinatorial} Schubert class.

To give the combinatorial description of the (combinatorial) Schubert class $\tau_u \in H_{\alpha}^*(\Gamma)$, we start by recalling some results concerning the combinatorics of the Weyl group $W$. Most results are valid for general Coxeter groups, and details can be found, for example, in \cite{BB}.

Let $R$ be a root system of rank $n$, let $W$ be its Weyl group,
and let $\mathcal{B} = \{ \alpha_1,\ldots,\alpha_n\}$ be a choice
of simple roots. Then $\mathcal{B}$ is a basis of $\ft^*$,
the dual of the Lie algebra of $T$. For a non-zero vector $\beta \in \ft^*$
we say that $\beta \succ 0$ if $\beta$ is in the non-negative cone
over $\{ \alpha_1,\ldots,\alpha_n\}$, in other words, if the coordinates of
$\beta$ in the basis $\mathcal{B}$ are non-negative.

For $i=1,\ldots,n$,
let $s_i =s_{\alpha_i} \in W$ be the reflection generated $\alpha_i$.
A word of length $m$ is an array $I = [i_1,i_2,\ldots,i_m]$ with
entries (letters) from $\{0,1,\ldots,n\}$. To each nonempty word
$I$ we associate the Weyl group element $s_I = s_{i_1}s_{i_2}\ldots
s_{i_m}$. (If the word $I$ is empty, then $s_I$ is the identity.) A
subword of a word $I$ is a word $J=[\epsilon_1 i_1, \epsilon_2
i_2,\ldots,\epsilon_m i_m]$, with $\epsilon_k \in \{0,1\}$ for all
$k=1,\ldots,m$. A word $I$ is \emph{reduced} if $s_{J} \neq
s_I$ for all subwords $J$ of $I$ other than $I$ itself. If $v \in W$,
then $v=s_{i_1}s_{i_2}\ldots s_{i_m}$ is called a \emph{decomposition} for $v$,
and $I = [i_1,i_2,\ldots,i_m]$ is a word for $v$. The decomposition is \emph{reduced}
if the word is reduced. All reduced words for $v \in W$ have the same number of letters,
and this common number, denoted by $\ell(v)$, is the \emph{length}
of $v$.

We define a partial order on $W$ as follows. If $v=us_{\beta}$, with
$\beta \succ 0$, we define $u \prec v$ if $u\beta \succ 0$.
The \emph{(strong) Bruhat order} on $W$ is the transitive closure of $\prec$. In other
words, $u \preccurlyeq v$ if and only if there exists an ascending chain
\begin{equation} \label{eq:relchain}
u = u_0 \xrightarrow{*s_{\beta_1}} u_1 \xrightarrow{*s_{\beta_2}} u_2
\longrightarrow \dotsb \xrightarrow{*s_{\beta_m}} u_m = v \; .
\end{equation}
where $p \xrightarrow{*s_{\beta}} q$ means $q=ps_{\beta}$. An
equivalent definition can be given in terms of words and subwords:
$u \preccurlyeq v$ if and only if every reduced word $I$ for $v$ has
a subword that is a word for $u$. Let $\mathcal{F}_u = \{ v \, | \, u \preccurlyeq v\}$ be the flow-up from $u$ under the strong Bruhat order.

For $u \in W$, we define
$$\Lambda_u^{-} \, = \,  \prod_{}^{} \{ \beta \; | \;  \beta \succ 0 \text{ and } u^{-1}\beta \prec 0\; \}$$

If $I = [i_1,\ldots,i_m]$ is a reduced word for $v$, then the positive roots that are sent by $v^{-1}$ into negative roots are
$$
\alpha_{i_1}\; , \;\; s_{i_1}\alpha_{i_2}\; , \;\;
s_{i_1}s_{i_2}\alpha_{i_3}\; , \;\; \ldots \; , \;\; s_{i_1}\dotsb
s_{i_{m-1}} \alpha_{i_m} \; ,
$$
hence
$$\Lambda_v^{-} = \prod_{j=1}^m s_{i_1}s_{i_2}\dotsb s_{i_j-1}
\alpha_{i_j} \; .$$

\begin{example}{\rm
  We illustrate the results above for $G=SL_{n}(\CC)$, the special linear group of complex matrices with determinant equal to one, with the Borel subgroup $B$ being the subgroup of upper triangular matrices. Then $M=G/B$ is the manifold of complete flags in $\CC^{n}$. The fixed points of the $T-$action are the classes of permutation matrices, hence $M^T$ is in a bijective correspondence with $W=S_{n}$, the set of permutations of $[n]=\{1,\ldots,n\}$. We use the one-line notation $u=u(1)u(2)\ldots u(n)$ for permutations. The simple roots are $\alpha_1=x_1\!-\!x_2$, \ldots, $\alpha_{n-1} = x_{n-1}\!-\!x_{n}$, the positive roots are $\alpha_{ij} = x_i\!-\!x_j = \alpha_i\! + \!\dotsb \!+\! \alpha_{j-1}$, for $1 \leqslant i < j \leqslant n$, and the Weyl group $W=S_n$ acts on roots by permuting the $x$ variables. The reflection $s_{\alpha_{ij}}$ corresponds to the transposition $(i,j)$ that swaps $i$ and $j$; if $j=i\!+\!1$, we denote $s_{\alpha_{ij}}$ by $s_i$ and call it a \emph{simple} transposition. By convention, $s_0$ is the identity of $S_n$.

  For a permutation $v \in S_n$, the set of positive roots that are sent by $v^{-1}$ into negative roots corresponds bijectively to the set of inversions-as-values,
$$\Inv(v)=\{ (v(j),v(i)) \; | \; i < j \; , \;  v(i)>v(j)\} \; . $$
If $I = [i_1,\ldots,i_m]$ is a reduced word for $v$, then
$$
  \Lambda_v^{-} =
   \prod_{(a,b) \in \; \Inv(v)} (x_a-x_b)
    =  \prod_{k=1}^m s_{i_1}\dotsb s_{i_{k-1}} \alpha_{i_k}
   \; .
$$

If $v=u(i,j)$ with $i<j$, then $u \prec v$ in the strong Bruhat order if and only if $u(i) < u(j)$. In other words, $u \prec v = u(i,j)$ if $(u(i),u(j)) \not\in \Inv(u)$, or, equivalently, if $(v(i),v(j)) \in \Inv(v)$.
  }
\end{example}

The description of the combinatorial Schubert class $\tau_u$ is the following.

\begin{prop}{\rm
For every $u \in W$, the class $\tau_u \in
H_{\alpha}^*(\Gamma)$ is the unique class satisfying the following
properties:
\begin{enumerate}
  \item \label{cond:1} $\tau_u$ is homogeneous, of degree $\ell(u)$, the length of $u$;
  \item \label{cond:2} $\tau_u$ is supported on $\mathcal{F}(u)$, the flow-up from $u$;
  \item \label{cond:3} $\tau_u$ is normalized by $\tau_u(u) =  \Lambda_u^{-}\; .$
\end{enumerate}}
\end{prop}

\section{Positive Formula Using Subwords}
\label{sec:subwordformula}

In this section we recall Billey's integral positive formula.
Let $u,v \in W$, with $u \preccurlyeq v$, and let $I =
[i_1,i_2,\ldots,i_m]$ be a reduced word for $v$. For a subword
$J=[\epsilon_1 i_1, \epsilon_2 i_2,\ldots,\epsilon_m i_m]$ of $I$,
define the \emph{subword contribution} $SC(J,I)$ as the following
product of positive roots:
\begin{align*}
  SC(J,I) = & \prod_{j=1}^m \left[ s_{i_1}s_{i_2}\dotsb s_{i_j-1} \alpha_{i_j}
  \right]^{\epsilon_j} = \prod_{ \substack{j=1 \\ \epsilon_j=1}}^m s_{i_1}s_{i_2}\dotsb
  s_{i_j-1} \alpha_{i_j} = \\
  = & \Lambda_v^{-} \cdot \prod_{\substack{j=1 \\ \epsilon_j
  =0}}^m \frac{1}{s_{i_1}s_{i_2}\dotsb s_{i_j-1}\alpha_{i_j}} \; ,
\end{align*}
hence $SC(J,I)$ is obtained from $\Lambda_v^{-}$ by canceling the positive roots generated by
the deleted letters ($\epsilon_j=0$). Let $\cR(u,I)$ be the set of
subwords of $I$ that, after deleting the zeroes, become
\emph{reduced} words for $u$.

\begin{theorem}[\cite{Bi}]{\rm Let $u,v \in W$ and let $I=[i_1,\ldots,i_m]$
be a reduced word for $v$. The value at $v$ of the Schubert class
$\tau_u$ is given by
\begin{equation}
\label{eq:posFormula}
  \tau_u(v) = \sum_{J \in \cR(u,I)} SC(J,I) = \sum_{J \in \cR(u,I)}
   \prod_{ \substack{j=1 \\ \epsilon_j=1}}^m s_{i_1}s_{i_2}\dotsb
  s_{i_j-1} \alpha_{i_j} \; .
\end{equation}}
\end{theorem}

\begin{example}
\label{ex:subwordContrib} {\em
  Let $u=2143$ and $v=3421$ in $S_4$.
   A reduced word for $v$ is
  $I=[2,1,3,2,3]$. There are two subwords of $I$ that are reduced
  words for $u$ after deleting all the zeroes:
  $J_1=[0,1,3,0,0]$ and $J_2=[0,1,0,0,3]$. Their contributions are
\begin{align*}
SC(J_1,I)=& [\alpha_{2}]^0 \cdot [s_2\alpha_1]^1 \cdot
[s_2s_1\alpha_3]^1 \cdot [s_2s_1s_3\alpha_2]^0 \cdot
[s_2s_1s_3s_2\alpha_3]^0 = \nonumber \\ = & [s_2\cdot (x_1-x_2)]
\cdot [s_2s_1\cdot (x_3-x_4)] = (x_1-x_3)(x_2-x_4) = \nonumber \\
= & (\alpha_1+\alpha_2)(\alpha_2+\alpha_3) \\
%
SC(J_2,I)=& [\alpha_{2}]^0 \cdot [s_2\alpha_1]^1 \cdot
[s_2s_1\alpha_3]^0 \cdot [s_2s_1s_3\alpha_2]^0 \cdot
[s_2s_1s_3s_2\alpha_3]^1 = \nonumber \\ = & [s_2\cdot (x_1-x_2)]
\cdot [s_2s_1s_3s_2\cdot (x_3-x_4)] = (x_1-x_3)(x_1-x_2) = \nonumber \\
= & (\alpha_1+\alpha_2)\alpha_1 \; .
\end{align*}
Therefore
\begin{equation*}
  \tau_{2143}(3421) = SC(J_1,I)+SC(J_2,I) =
  (\alpha_1+\alpha_2)(\alpha_1+\alpha_2+\alpha_3) \; .
\end{equation*}
  }
\end{example}

\section{The Goldin-Tolman Formula}
\label{sec:GTformula}

Goldin and Tolman (\cite{GT}) have recently announced a formula for computing the values of what they call \emph{canonical classes}, valid for more general spaces\footnote{Allen Knutson has informed me that the flag manifold case of this formula also appeared, in implicit form, in his earlier paper \cite{Kn}.}. For flag manifolds, their classes satisfy conditions~\ref{cond:1} and \ref{cond:3} above, but condition~\ref{cond:2} is replaced by $\tau_{u}(v) = 0 $ if $\ell(v) \leqslant \ell(u)$ and $v \neq u$. However, since the length function is strictly increasing with respect to the Bruhat order (that is, if $u \prec v$ then $\ell(u) < \ell(v)$), it turns out that in this case the canonical generators are the equivariant Schubert classes, and their formula can be used to compute the values $\tau_u(v)$ for $u,v \in W$.

The Goldin-Tolman formula involves two more ingredients.
The first ingredient is the subgraph $\Gamma_0$ of $\Gamma$, having
the same vertices as $\Gamma$, but only the edges $e = (u,v)$ of
$\Gamma$ for which $\ell(v) = \ell(u) \pm 1$. If $u$ and $v$ are
elements in $W$, let $\Sigma(u,v)$ be the set of ascending
chains in $\Gamma_0$ from $u$ to $v$. These chains are the maximal
length ascending chains in $\Gamma$ from $u$ to $v$. Every edge of
such a chain will be considered oriented, with the orientation that
makes it an ascending edge.

The second ingredient is a moment map. The Lie algebra $\fg$ of $G$ can be canonically identified with its dual $\fg^*$ using the Killing form, and that allows us to regard $\ft^*$ as a subspace of $\fg^*$. Let $\eta \in \ft^* \subset \fg^*$ be in the interior of the positive Weyl chamber, and let $\cO_{\eta} = G \cdot \eta \subset \fg^*$ be the co-adjoint orbit through $\eta$. Then the stabilizer of $\eta$ is $B$, and hence $\cO_{\eta} \simeq G/B=M$ as $T-$spaces.

The co-adjoint orbit $\cO_{\eta}$ is a Hamiltonian $G-$space, with moment map given by the inclusion $\cO_{\eta} \hookrightarrow \fg^*$. Therefore it is also a Hamiltonian $T-$space, with moment map given by inclusion followed by projection onto $\ft^*$. With the symplectic structure induced by the identification $G/B \simeq \cO_{\eta}$, the flag manifold $M=G/B$ is a Hamiltonian $T-$space, with moment map $\phi_{\eta} \colon G/B \to \ft^*$. If $P_{w}$ is the fixed point corresponding to the element $w\in W$ of the Weyl group, then $\phi_{\eta}(P_w) = w\cdot \eta$. Identifying the fixed point $P_w$ with the Weyl group element $w$, we get a map $\phi_{\eta} \colon W \to \ft^*$, given by $\phi_{\eta}(w) = w\eta$.

Applying the Goldin-Tolman formula (\cite{GT}), we get the following result.

\begin{thm}{\rm If $\eta \in \ft^*$ is in the positive Weyl chamber, then
\begin{equation}\label{eq:GTformula}
  \tau_u(v) =  \sum_{\gamma \in
  \Sigma(u,v)} E_{\eta}(\gamma) \; ,
\end{equation}
where, for every ascending chain of maximal length
\begin{equation}\label{eq:maxChain}
  \gamma : \; u=u_0 \longrightarrow u_1 \longrightarrow
   \dotsb \longrightarrow u_m = v \;
\end{equation}
in $\Sigma(u,v)$, the contribution $E_{\eta}(\gamma)$ is given by
\begin{equation}
\label{eq:genElambda}
  E_{\eta}(\gamma) = \Lambda_v^{-} \prod_{k=1}^m
 \Bigl(  \frac{1}{\alpha(u_{k-1},u_k)}\cdot \frac{\phi_{\eta}(u_k)-\phi_{\eta}(u_{k-1})}
 {\phi_{\eta}(v) - \phi_{\eta}(u_{k-1})} \Bigr) \; .
\end{equation}}
\end{thm}

There are two important features in this formula. The first is that
each term depends on $\eta$, but the sum doesn't. The second
is that each term is a rational expression, but the sum is a
polynomial. Therefore a lot of cancelations must occur when summing,
and a first indication of how that happens is the following.

\begin{lemma}\label{lem:1}{\rm
  If $p,q \in W$ and $q = ps_{\beta}$ with $\beta \succ 0$, then
  \begin{equation}
    \frac{\phi_{\eta}(q)-\phi_{\eta}(p)}{\alpha(p,q)} <0\; .
  \end{equation}}
\end{lemma}

\begin{proof}
$$\frac{\phi_{\eta}(q)-\phi_{\eta}(p)}{\alpha(p,q)} =   \frac{ps_{\beta}\eta - p\eta}{p\beta} = - \langle \eta,\beta \rangle  = -\frac{2(\eta,\beta)}{(\beta,\beta)} < 0 \; ,$$
since $\beta$ is a positive root and $\eta$ is in the positive Weyl chamber.
\end{proof}

Hence, if $\gamma$ is the maximal chain \eqref{eq:maxChain} and $u_k = u_{k-1}s_{\beta_k}$ with $\beta_k \succ 0$, then
\begin{equation*}
E_{\eta}(\gamma) = \Lambda_v^{-} \prod_{k=1}^m \frac{\langle \eta,\beta_k \rangle}{u_{k-1}\eta -v\eta} \; .
\end{equation*}

\begin{lemma}\label{lem:2}{\rm
   Let $p \prec q$, let $\eta$ be in the closure of the positive Weyl chamber and let $\beta = p\eta -q\eta$. Then
   $\beta \succcurlyeq 0$. If $\beta \succ 0$, then $p^{-1}\beta \succ 0$ and $q^{-1}\beta \prec 0$.
   }
\end{lemma}

\begin{proof}
 If $p \prec q$, then there exists an ascending chain
 $$p = p_0 \xrightarrow{*s_{\beta_1}} p_1 \xrightarrow{*s_{\beta_2}} p_2 \to \dotsb \xrightarrow{*s_{\beta_m}} p_m = q$$
 with $\beta_k \succ 0$ and $p_{k-1} \beta_k \succ 0$. Then
\begin{align}\label{eq:nonnegcone}
  p\eta -q\eta = & (p \eta - p_1 \eta) + \dotsb + (p_{k-1}\eta -p_k\eta) + \dotsb + (p_{m-1}\eta -q\eta) = \\
               = & \langle \eta,\beta_1\rangle p_{0}\beta_1 + \dotsb + \langle \eta,\beta_k\rangle p_{k-1}\beta_k + \dotsb + \langle \eta,\beta_m\rangle p_{m-1}\beta_m \nonumber \; ,
\end{align}
and each of the terms in the last sum is in the non-negative cone generated by the simple roots, hence $\beta \succcurlyeq 0$. Moreover, $\beta = 0$ if and only if all the factors $\langle \eta, \beta_k \rangle$ are zero.

We have $p^{-1}\beta = \eta - p^{-1}q\,\eta = id\; \eta - p^{-1}q\,\eta$,
where $id$ is the identity in the Weyl group. Since $id \prec p^{-1}q$, if $\beta \succ 0$ then $p^{-1}\beta \succ 0$, and similarly, $q^{-1}\eta = - (id\; \eta - q^{-1}p\,\eta) \prec 0$.
\end{proof}

Therefore $E_{\eta}(\gamma)$
is a ratio of two polynomials in the simple roots $\alpha_1,\ldots,\alpha_n$, homogeneous and with non-negative coefficients. Moreover, if $u_{k-1} \eta - v\eta$ is a multiple of a root, then it cancels (over $\RR_{>0}$) one of the factors of $\Lambda_v^{-}$.

Let $\{ \omega_1,\ldots,\omega_n\}$ be the fundamental weights corresponding to the simple roots $\{\alpha_1,\ldots,\alpha_n\}$, defined by the conditions
$$\langle \omega_i,\alpha_j \rangle = \delta_{ij} \quad \text{ or, equivalently, } \quad (\omega_i,\alpha_j) = \frac{1}{2}(\alpha_i,\alpha_i) \delta_{ij} \; .$$
Then $\{ \omega_1,\ldots,\omega_n\}$ is a basis of $\ft^*$, and
$\eta = \mu_1 \omega_1 + \dotsb + \mu_n \omega_n \in \ft^*$ is in the positive Weyl chamber if and only if $\mu_j>0$ for all $j=1,\ldots,n$.

Then the contribution of the chain $\gamma$ is
\begin{equation}\label{eq:chainContrib_mu}
  E_{\mu}(\gamma) = E_{\mu_1\omega_1+ \dotsb + \mu_n\omega_n}(\gamma) = \Lambda_v^{-} \prod_{k=1}^m
\frac{\sum_{i=1}^{n} \langle \omega_i,\beta_k\rangle \mu_i}
{\sum_{i=1}^n \left( \sum_{j=k}^m \langle \omega_i,\beta_j \rangle u_{j-1}\beta_j \right) \mu_i} \; .
\end{equation}

\begin{example}\label{exm:simpleA2}
  {\rm
  Let $M=SL_3(\CC)/B$ be the manifold of complete flags in $\CC^3$, let $u=s_1 = 213$ and $v=s_1s_2s_1= 321$. Then $\ell(321) = 3$, and
$$
  \Lambda_v^{-} =  \alpha_1\alpha_2(\alpha_1+\alpha_2) \; .
$$
There are two ascending chains of maximal length from $u$ to $v$,
\begin{align*}
  \gamma_1 \; \colon &  \quad s_1 \xrightarrow{*s_{\alpha_1+\alpha_2}} s_2s_1
  \xrightarrow{*s_{\alpha_2}} s_1s_2s_1 \\
  \gamma_2 \; \colon &  \quad s_1 \xrightarrow{*s_{\alpha_2}} s_1s_2
  \xrightarrow{*s_{\alpha_1}} s_1s_2s_1\; ,
\end{align*}
and the contributions of the two chain are 
\begin{align*}
  E_{\mu}(\gamma_1) = & \alpha_1\alpha_2(\alpha_1\!+\!\alpha_2) \! \cdot \! \frac{\mu_1\!+\!\mu_2}{\alpha_2\mu_1\!+\!(\alpha_1\!+\!\alpha_2)\mu_2} \! \cdot \! \frac{\mu_2}{\alpha_1\mu_2} =
   \frac{\alpha_2(\alpha_1\!+\!\alpha_2)(\mu_1\!+\!\mu_2)}{\alpha_2\mu_1\!+\!(\alpha_1\!+\!\alpha_2)\mu_2}
  \\
%
  E_{\mu}(\gamma_2) = & \alpha_1\alpha_2(\alpha_1\!+\!\alpha_2) \! \cdot \! \frac{\mu_2}{\alpha_2\mu_1\!+\!(\alpha_1\!+\!\alpha_2)\mu_2} \! \cdot \! \frac{\mu_1}{\alpha_2\mu_1} =  \frac{\alpha_1(\alpha_1\!+\!\alpha_2)\mu_2}{\alpha_2\mu_1\!+\!(\alpha_1\!+\!\alpha_2)\mu_2} \; .
\end{align*}
Then
\begin{equation*}
  \tau_{u}(v) = E_{\mu}(\gamma_1) + E_{\mu}(\gamma_2) = \alpha_1 +
  \alpha_2 \; .
\end{equation*}
  }
\end{example}

\section{Limits and Chain Contributions}
\label{sec:limits}

Note that, in Example~\ref{exm:simpleA2}, both $E_{\mu}(\gamma_1)$ and $E_{\mu}(\gamma_2)$
are rational expressions in $\alpha$'s and $\mu$'s, but their sum is a
polynomial expression in the $\alpha$ variables only. In this section we show that we can
eliminate the $\mu$ variables in \eqref{eq:chainContrib_mu} by sending them to 0, one component
at a time.

The contribution \eqref{eq:genElambda} of the chain \eqref{eq:maxChain} can be written as
$$E_{\mu}(\gamma) = \Lambda_v^{-} \prod_{k=1}^{m} \frac{Q_{\mu}(u_{k-1},u_k,v)}{\alpha(u_{k-1},u_k)} \; ,$$
where, for $p \prec q \preccurlyeq r$ in $W$, we define
\begin{equation}
  Q_{\mu}(p,q,r) = \frac{p\, \eta - q \eta}{p\, \eta-r \eta} = \frac{\sum_{i=1}^n (p\, \omega_i -q \omega_i)\mu_i}{\sum_{i=1}^n (p\, \omega_i -r \omega_i)\mu_i}\; .
\end{equation}
\begin{lemma}\label{lem:4}{\rm
  For $p,q \in W$, let
  $$h(p,q) = \min \{ i \; | \; p\,\omega_i \neq q\omega_i\} \; ,$$
  with the convention that $h(p,p)=\infty$.
  \begin{enumerate}
    \item If $p^{-1}q = s_{i_1}s_{i_2}\dotsb s_{i_k}$ is a reduced decomposition, then
    $$h(p,q) = \min\{i_1,\ldots,i_k\} \; .$$
    \item If $p \prec q \prec r$ and $p\, \omega_i = r\omega_i$, then $p\,\omega_i = q\omega_i = r \omega_i$.
    \item If $p \prec q \prec r$, then
     \begin{equation*}
       h(p,r) \leqslant h(p,q) \quad \text{ and } \quad h(q,r) \leqslant h(p,q) \; .
     \end{equation*}
     \item If $q = ps_{\beta}$ with $\beta \succ 0$ and $p\beta \succ 0$, then $h(p,q) = h(\beta)$, where, for nonzero $\beta = \beta^1 \alpha_1 + \dotsb + \beta^n \alpha_n$ in $\ft^*$, we define $h(\beta) = \min\{ i \; | \; \beta^i \neq 0\}$.
  \end{enumerate}}
\end{lemma}

\begin{proof}
Let $j=h(p,q)$ and $\omega_i$ a fundamental weight. Then \eqref{eq:nonnegcone} becomes
\begin{equation}\label{eq:minleter}
  \omega_i - p^{-1}q \omega_i = \delta_{i_1,\,i\,} \alpha_{i_1} + \delta_{i_2,\,i\,}s_{i_1}\alpha_{i_2} + \dotsb + \delta_{i_k,\,i\,} s_{i_1}s_{i_2}\dotsb s_{i_{k-1}}\alpha_{i_k} \; .
\end{equation}
If $i<j$, then $\omega_i = p^{-1}q \omega_i$, and since all the nonzero terms in \eqref{eq:minleter} are positive roots, it follows that $s_i$ can not appear in the reduced decomposition of $p^{-1}q$. If $i=j$, then $\omega_j \neq p^{-1}q \omega_j$, hence $s_j$ must appear in the decomposition of $p^{-1}q$. Therefore $j=h(p,q)$ is the minimal letter that appears in a reduced word for $p^{-1}q$ and $p^{-1}q$ is in the subgroup of $W$ generated by $s_j, s_{j+1},\ldots,s_n$.

  For every $i=1,\ldots,n$ we have $p\, \omega_i -r \omega_i = (p\,\omega_i - q \omega_i) + (q \omega_i - r\omega_i)$, and since $p \prec q \prec r$, both $p\,\omega_i - q \omega_i$ and $q \omega_i - r\omega_i$ are in the non-negative cone generated by the simple roots. Therefore, if $p\, \omega_i -r \omega_i = 0$, then $p\,\omega_i - q \omega_i=0$ and $q \omega_i - r\omega_i=0$. Then $\{ i \; | \; p\, \omega_i = r \omega_i\} \subset \{ i \; | \; p\, \omega_i = q \omega_i\}$, which implies $h(p,r) \leqslant h(p,q)$. Similarly we get $h(p,r) \leqslant h(q,r)$.

  If $q = ps_{\beta}$, then
$$ \{ i \; | \; p\,\omega_i \neq ps_{\beta}\omega_i\} =
   \{ i \; | \; \omega_i \neq s_{\beta}\omega_i\} =
   \{ i \; | \; (\omega_i, \beta) \neq 0 \} =
   \{ i \; | \; \beta^i \neq 0 \}\;, $$
hence $h(p,q) = h(\beta)$.
\end{proof}

Consequently, if the coefficient of some $\mu_i$ in the denominator is zero, then the coefficient of $\mu_i$ in the numerator is also zero. This implies that
$$Q_{\mu}^0(p,q,r) \stackrel{\text{def}}{==} \lim_{\mu_1 \to 0} \Bigl(\, \lim_{\mu_2 \to 0} \Bigl(\,
  \dotsb \Bigl(\, \lim_{\mu_n \to 0} Q_{\mu}(p,q,r) \,\Bigr) \dotsb
  \,\Bigr) \,\Bigr)$$
is well defined. Moreover,
$$Q_{\mu}^0(p,q,r) =
\left\{
\begin{aligned}
  \frac{p\,\omega_i - q\omega_i}{p\,\omega_i - r\omega_i} \; , \; \text{ if } h(p,r) = h(p,q)=i \\
  0 \; , \; \text{ if } h(p,r) < h(p,q) \; .
\end{aligned}
\right.
$$
In particular, let $p=u_{k-1}$, $q=u_k=u_{k-1}\beta_k$, and $r=v$. Then
$$\frac{Q_{\mu}^0(u_{k-1},u_k,v)}{\alpha(u_{k-1},u_k)} =
\left\{
\begin{aligned}
  \frac{\langle \omega_i , \beta_k \rangle}{u_{k-1}\omega_i - v\omega_i} \; , \; \text{ if }
  h(u_{k-1},v) = h(u_{k-1},u_k)=i \\
  0 \; , \; \text{ if } h(u_{k-1},v) < h(u_{k-1},u_k) \; .
\end{aligned}
\right.
$$
Therefore
$$E(\gamma) \stackrel{\text{def}}{==} \lim_{\mu_1 \to 0} \Bigl(\, \lim_{\mu_2 \to 0} \Bigl(\,
  \dotsb \Bigl(\, \lim_{\mu_n \to 0} E_{\mu}(\gamma) \,\Bigr) \dotsb
  \,\Bigr) \,\Bigr)$$
is well-defined, and is zero whenever $h(u_{k-1},v) < h(u_{k-1},u_k)$ for some $k$. If $h(u_{k-1},u_k) = h(u_{k-1},v)=i_k$ for all $k=1,\ldots,m$, then
\begin{equation*}
  E(\gamma) = \Lambda_v^{-} \prod_{k=1}^m
\frac{\langle \omega_{i_k}, \beta_k \rangle}{u_{k-1}\omega_{i_k} - v\omega_{i_k}} \neq 0 \; .
\end{equation*}

We now take a closer look at the set of chains
$$\cC_0(u,v) = \{ \gamma \in \Sigma(u,v) \; | \; h(u_{k-1},u_k) = h(u_{k-1},v) \text{ for all } k=1,\ldots,m\} $$
with non-zero contribution after taking the limits. These chains have a simple description.

\begin{lemma}
  {\rm Let  $\gamma \in \Sigma(u,v)$ be the ascending chain of maximal length
  $$\gamma : \; u=u_0 \xrightarrow{*s_{\beta_1}} u_1 \xrightarrow{*s_{\beta_2}}
   \dotsb \xrightarrow{*s_{\beta_m}} u_m = v \; . $$
 Then $\gamma \in \cC_0(u,v)$ if and only if
 $h(\beta_1) \leqslant h(\beta_2) \leqslant \dotsb \leqslant h(\beta_m)$.
}
\end{lemma}

\begin{proof}

If $\gamma \in \cC_0(u,v)$, then, using Lemma~\ref{lem:4}, we get
$$h(\beta_{k+1}) = h(u_k,u_{k+1}) = h(u_k,v) \geqslant h(u_{k-1},v) = h(u_{k-1},u_k) = h(\beta_k)\; $$
for all $k=1,\ldots,m-1$, proving one implication.

Conversely, let $\gamma \in \Sigma(u,v)$ such that $h(\beta_1) \leqslant h(\beta_2) \leqslant \dotsb \leqslant h(\beta_m)$. Let $k \in [m]$ be fixed.  If $i < h(\beta_k)$, then $i < h(\beta_j) = h(u_{j-1},u_j)$ for all $j \geqslant k$, hence $u_{j-1}\, \omega_i = u_j \,\omega_i$ for all $j \geqslant k$. Therefore $u_{k-1}\, \omega_i = u_m \, \omega_i = v\omega_i$, hence $h(u_{k-1},v) \geqslant h(\beta_k) = h(u_{k-1},u_k)$. But $h(u_{k-1},v) \leqslant h(u_{k-1},u_k)$ by Lemma~\ref{lem:4}, and the double inequality implies that $h(u_{k-1},v) = h(u_{k-1},u_k)$ for all $k \in [m]$, and therefore $\gamma \in \cC_0(u,v)$.
\end{proof}

\begin{exm*}{\rm
  Returning to Example~\ref{exm:simpleA2}, we see that
\begin{align*}
  E(\gamma_1) = & \lim_{\mu_1 \to 0} \Bigl(\, \lim_{\mu_2 \to 0} \frac{\alpha_2(\alpha_1\!+\!\alpha_2)(\mu_1\!+\!\mu_2)}{\alpha_2\mu_1\!+\!(\alpha_1\!+\!\alpha_2)\mu_2} \, \Bigr) = \alpha_1 + \alpha_2 \neq 0\\
  E(\gamma_2) = & \lim_{\mu_1 \to 0} \Bigl(\, \lim_{\mu_2 \to 0}
  \frac{\alpha_1(\alpha_1\!+\!\alpha_2)\mu_2}{\alpha_2\mu_1\!+\!(\alpha_1\!+\!\alpha_2)\mu_2}
  \Bigr) = 0\; ,
\end{align*}
as expected, since, for $\gamma_1$ we have $h(\alpha_1\! + \!\alpha_2) = 1 \leqslant 2 = h(\alpha_2)$ and for $\gamma_2$ we have $h(\alpha_2) = 2 \nleqslant 1 = h(\alpha_1)$.}
\end{exm*}

\section{Positive Formula Using Chains}
\label{sec:chainformula}

Up to this point, there have been no restrictions regarding the order in which the simple roots $\alpha_1,\ldots,\alpha_n$ are listed. The next results, however, requires that the simple roots be listed in a specific order, as shown in the following Dynkin diagrams:
\begin{equation}\label{eq:Dynkin}
  \begin{aligned}
  \text{Type $A$: } &  \quad \alpha_1 -\!\!\!-\!\!\!-\!\!\!- \alpha_2 -\!\!\!-\!\!\!-\!\!\!-
   \dotsb -\!\!\!-\!\!\!-\!\!\!- \alpha_{n-1} -\!\!\!-\!\!\!-\!\!\!- \alpha_n \\
  \text{Type $B$: } &  \quad \alpha_1 -\!\!\!-\!\!\!-\!\!\!-
   \alpha_2 -\!\!\!-\!\!\!-\!\!\!- \dotsb -\!\!\!-\!\!\!-\!\!\!- \alpha_{n-1} =\!=\!=\!\!\!\!\!\!\!\!> \;\; \alpha_n \\
  \text{Type $C$: } &  \quad \alpha_1 -\!\!\!-\!\!\!-\!\!\!-
   \alpha_2 -\!\!\!-\!\!\!-\!\!\!- \dotsb -\!\!\!-\!\!\!-\!\!\!- \alpha_{n-1} =\!=\!=\!\!\!\!\!\!\!\!< \;\; \alpha_n
\end{aligned}
\end{equation}

We are now ready to formulate the main result of this paper. Let $\cC_0(u,v)$ be
the set of ascending chains of maximal length
  \begin{equation}\label{eq:maxacsendchain}
    \gamma \colon \; u=u_0 \xrightarrow{*s_{\beta_1}} u_1 \xrightarrow{*s_{\beta_2}} u_2 \to \dotsb \to u_{m-1} \xrightarrow{*s_{\beta_m}} u_m = v
  \end{equation}
  that satisfy the condition
  \begin{equation}\label{eq:convexchains}
    h(\beta_1) \leqslant h(\beta_2) \leqslant \dotsb \leqslant h(\beta_m) \; ,
  \end{equation}
  and for each such a chain $\gamma$, let
  \begin{equation}\label{eq:Emu0}
    E(\gamma) = \Lambda_v^{-} \prod_{k=1}^m
\frac{\langle \omega_{i_k}, \beta_k \rangle}{u_{k-1}\omega_{i_k} - v\omega_{i_k}} \; ,
  \end{equation}
  where $i_k = h(\beta_k)$ for $k=1,\ldots,m$.

\begin{thm}\label{thm:main-general}{\rm If $u, v \in W$ then:

\textbf{[Part I]} The restriction of $\tau_u$ at $v$ is given by
  \begin{equation}\label{eq:genpositive}
    \tau_u(v) = \sum_{\gamma \in\, \cC_0(u,v)} E(\gamma) =
    \sum_{\gamma \in\, \cC_0(u,v)} \Lambda_v^{-} \prod_{k=1}^m
\frac{\langle \omega_{i_k}, \beta_k \rangle}{u_{k-1}\omega_{i_k} - v\omega_{i_k}} \; .
  \end{equation}

\textbf{[Part II]} If the simple roots are ordered as in \eqref{eq:Dynkin}, then:
  \begin{description}
  \item[In types $A$ and $C$]
  $$E(\gamma) \in \ZZ_{\geqslant 0}[\alpha_1,\ldots,\alpha_n] \; ,$$
  hence \eqref{eq:genpositive} is a positive integral formula for computing $\tau_u(v)$.
  \item[In type $B$]
  $$E(\gamma) \in \frac{1}{2^m}\ZZ_{\geqslant 0}[\alpha_1,\ldots,\alpha_n] \; ,$$
  hence \eqref{eq:genpositive} is a positive, but only rational, formula for $\tau_u(v)$.
  \end{description}

\textbf{[Part III]} If the simple roots are ordered as in \eqref{eq:Dynkin}, then:

  \begin{description}
    \item[In type $A$] For every $u$ and $v$, there is a particular choice of a reduced word $I$ for $v$, and a bijection $F_I \colon \cC_0(u,v) \to \cR(u,I)$ such that $E(\gamma) = SC(F_I(\gamma),I)$. In other words, the positive formula using chains is an alternative formulation of Billey's positive formula using subwords.
        \item[In type $C$] There exist $u$ and $v$ such that, whichever reduced word $I$ for $v$ we choose, the number of subwords in $\cR(u,I)$ is not the same as the number of chains in $\cC_0(u,v)$. Therefore, in type $C$, the two positive formulas are not manifestly equivalent.
  \end{description}
  }
\end{thm}

\begin{proof}
  The right hand side of the Goldin-Tolman formula \eqref{eq:GTformula} is independent of the particular choice of $\eta$ in the positive chamber (and hence of $\mu$ with strictly positive components). Each term in the right hand side is well-defined after sending the components of $\mu$ to zero, one at a time in reverse order. If $\gamma \not\in \cC_0(u,v)$, then the limit is zero, and if $\gamma \in \cC_0(u,v)$, then the limit is $E(\gamma)$, and this proves the statement in Part I.

  With the ordering \eqref{eq:Dynkin}, if $\beta$ is a positive root and $h(\beta)=i$, then $\langle \omega_i, \beta \rangle$ is a positive integer, so the product in the numerator of \eqref{eq:Emu0} is a positive integer. To prove the statements of Part II, we show that:
  \begin{enumerate}
    \item Each factor in the denominator of \eqref{eq:Emu0} cancels (over $\ZZ_{>0}$ in types $A$ and $C$, and over $\frac{1}{2}\ZZ_{>0}$ in type $B$) one of the factors in $\Lambda_{v}^{-}$, and
    \item The factors in the denominator of \eqref{eq:Emu0} are mutually independent.
  \end{enumerate}
  These statements are consequences of the following technical lemma, which can be verified through straightforward computations.
  \begin{lemma}\label{lem:technlemma}{\rm
    If $u \in W$ and $j=h(id,u)$, then:
  \begin{description}
    \item[In type $A$]{} $u$ has a reduced decomposition of the form
    \begin{equation}\label{eq:decompA}
      u=s_ks_{k-1}\dotsb s_jw
    \end{equation}
    with $j \leqslant k \leqslant n$ and $h(id,w)>j$.
    \smallskip
    \item[In type $B$ or $C$] $u$ has either a reduced decomposition \eqref{eq:decompA}
    or a reduced decomposition of the form
    \begin{equation}\label{eq:decompBC}
      u = s_t s_{t+1} \dotsb s_{n-1}s_ns_{n-1} \dotsb s_{j+1}s_j w \; ,
    \end{equation}
    with $j \leqslant t < n$ and $h(id,w)>j$, but not both.
  \end{description}

  If $u$ has a decomposition of the form \eqref{eq:decompA}, then
  $$\omega_j -u\omega_j = \left\{
  \begin{aligned}
    \alpha_j+ \dotsb + \alpha_k, & \text{ in types $A, C$ for $k\leqslant n$}  \\
                                 & \text{ and type $B$ if $k<n$} \\
    \alpha_j+ \dotsb + \alpha_{n-1}+ 2\alpha_n, & \text{ in type $B$ for $k=n$}
  \end{aligned}
  \right.$$

  If $u$ has a decomposition of the form \eqref{eq:decompBC} with $t>j$, then
  $$\omega_j -u\omega_j = \left\{
  \begin{aligned}
    \alpha_j+ \dotsb + \alpha_{t-1}+2(\alpha_t + \dotsb + \alpha_n), & \text{ in type $B$} \\
    \alpha_j+ \dotsb + \alpha_{t-1}+2(\alpha_t + \dotsb + \alpha_{n-1})+\alpha_n, & \text{ in type $C$}
  \end{aligned}
  \right.$$

  If $u$ has a decomposition of the form \eqref{eq:decompBC} with $t=j$, then
    $$
    \omega_j -u\omega_j =
    \left\{
    \begin{aligned}
      2(\alpha_j+\dotsb + \alpha_n) & \quad \text{, in type }B \\
      2(\alpha_j + \dotsb + \alpha_{n-1}) + \alpha_n & \quad \text{, in type } C
    \end{aligned}
    \right.
    $$
   If $u$ is not of the form
  $s_j\dotsb s_n \dotsb s_j w$ in type $B$, then $\omega_j -u\omega_j$ is a root.
In the excepted case, $\omega_j -u\omega_j$  is twice a root. In all cases $h(\omega_j - u \omega_j) = j$, and no distinct factors are multiples of each other.
}
  \end{lemma}

  The above lemma implies that
  $$u_{k-1}\omega_{i_k} -v \omega_{i_k} = u_{k-1}(\omega_{i_k} - u_{k-1}^{-1}v\omega_{i_k})$$
  is either a root (in types $A$ and $C$) or, possibly, twice a root (in type $B$). From Lemma~\ref{lem:2} we conclude that each factor of \eqref{eq:Emu0} cancels (over $\ZZ_{>0}$ in types $A$ and $C$, and over $\frac{1}{2}\ZZ_{>0}$ in type $B$) one of the factors in $\Lambda_{v}^{-}$.

  To finish the proof of Part II, we need to show that the factors are mutually independent. Suppose that is not true, and that for some chain $\gamma \in \cC_0(u,v)$ and some $k < t$, the vectors $u_{k-1}\omega_{i_k}-v \omega_{i_k}$ and $u_{t-1}\omega_{i_t}-v \omega_{i_t}$ are dependent. Then $\omega_{i_k} - v^{-1}u_{k-1}\omega_{i_k}$ and $\omega_{i_t} - v^{-1}u_{t-1}\omega_{i_t}$ are also dependent. Using Lemma~\ref{lem:technlemma} we get that
  $$i_k = h(\omega_{i_k} - v^{-1}u_{k-1}\omega_{i_k}) = h(\omega_{i_t} - v^{-1}u_{t-1}\omega_{i_t}) = i_t \; ,$$
  hence $i_k=i_t$. Let $j =i_k=i_t$ be the common value. Since $\omega_{j} - v^{-1}u_{k-1}\omega_{j}$ and $\omega_{j} - v^{-1}u_{t-1}\omega_{j}$ are dependent, Lemma~\ref{lem:technlemma} implies that they must be equal, and therefore
  $u_{k-1}\omega_j - u_{t-1}\omega_j = 0$. Then the proof of Lemma~\ref{lem:2} implies that $\langle \omega_j, \beta_k \rangle = 0$, and that is a contradiction, since $j=i_k = h(\beta_k)$.
 Therefore the factors in the denominator of \eqref{eq:Emu0} are independent, and that finishes the proof of Part II.

If $u=s_2$ and $v=s_1s_2s_1$ in $B_2$, then the chain
$$\gamma \, \colon \; s_2 \xrightarrow{*s_{\alpha_1+2\alpha_2}} s_1s_2 \xrightarrow{*s_{\alpha_1+\alpha_2}} s_1s_2s_1$$
is in $\cC_0(u,v)$, but
$$E(\gamma) = \alpha_1 (\alpha_1+\alpha_2)(\alpha_1+2\alpha_2) \cdot \frac{1}{2(\alpha_1+\alpha_2)}
\cdot \frac{1}{\alpha_1+2\alpha_2} = \frac{1}{2} \,\alpha_1\; ,$$
hence in type $B$ the formula \eqref{eq:genpositive} is not necessarily integral, and the denominators can be higher powers of 2.

The following example proves the Part III statement in type $C$.

Let $u=s_1$ and $v=s_1s_2s_1s_2$ in $C_2$. Then $I_1=[1,2,1,2]$ and $I_2=[2,1,2,1]$ are the only reduced words for $v$, and each of them has two subwords that are reduced words for $u$. Hence whether $I=I_1$ or $I=I_2$, there are two subwords in $\cR(u,I)$, hence two  terms in \eqref{eq:posFormula}, the positive formula using subwords.

There are four ascending chains of maximal length from $u$ to $v$:
\begin{equation}\label{eq:exampleC2}
\begin{aligned}
  \gamma_1 : & \quad s_1 \xrightarrow{*s_{\alpha_1+\alpha_2}} s_2s_1
                         \xrightarrow{*s_{\alpha_1+2\alpha_2}} s_1s_2s_1
                         \xrightarrow{*s_{\alpha_2}} s_1s_2s_1s_2\\
  \gamma_2 : & \quad s_1 \xrightarrow{*s_{\alpha_1+\alpha_2}} s_2s_1
                         \xrightarrow{*s_{\alpha_2}} s_2s_1s_2
                         \xrightarrow{*s_{\alpha_1}} s_1s_2s_1s_2\\
   \gamma_3 : & \quad s_1 \xrightarrow{*s_{\alpha_2}} s_1s_2
                         \xrightarrow{*s_{\alpha_1+\alpha_2}} s_2s_1s_2
                         \xrightarrow{*s_{\alpha_1}} s_1s_2s_1s_2\\
  \gamma_4 : & \quad s_1 \xrightarrow{*s_{\alpha_2}} s_1s_2
                         \xrightarrow{*s_{\alpha_1}} s_1s_2s_1
                         \xrightarrow{*s_{\alpha_2}} s_1s_2s_1s_2
\end{aligned}
\end{equation}
but only $\gamma_1$ satisfies \eqref{eq:convexchains}, hence $\cC_0(u,v) = \{ \gamma_1\}$ has only one chain.

All that remains to be proved is the Part III statement in type $A$. In Section~\ref{sec:TypeAchains} we give an explicit version of the positive formula \eqref{eq:posFormulaChains} in type $A$, in Section~\ref{sec:chains-subwords} we construct a general map connecting chains and subwords, and we finish the proof in Section~\ref{sec:equivalence}.
\end{proof}

The main result of Lemma~\ref{lem:technlemma} is that if $h(id,u)=j$, then $\omega_j - u\omega_j$ is a root (in type $A$ and $C$) or a multiple of a root (it type $B$, and, as one can check, for $G_2$). However, in type $D$ or for $F_4$, $\omega_j - u\omega_j$ is a sum of positive roots that is not a multiple of a root, and therefore, the corresponding chain contribution \eqref{eq:Emu0} is positive but not even polynomial.

\section{Type A Positive Formula Using Chains}
\label{sec:TypeAchains}

We apply the general result of Theorem~\ref{thm:main-general} to  $M=SL_{n}(\CC)/B$.

If $h(\beta)=k$, then $\beta=\alpha_{kj}$ for some $j \geqslant k$ and $\langle \omega_k, \beta \rangle = 1$. Hence all the factors in the numerator of \eqref{eq:Emu0} are equal to 1. The denominators can be computed using the following lemma.

\begin{lemma}{\rm
  Let $u \in S_n$. If $h(id, u)=j$, then
  \begin{equation}
    \omega_j - u \omega_j = x_j -x_{u(j)} \; .
  \end{equation}}
\end{lemma}

\begin{proof}
  Let $u = s_ks_{k-1}\dotsb s_jw$ be the decomposition \eqref{eq:decompA}. Then
  $$\omega_j - u \omega_j = \alpha_j + \dotsb \alpha_k = x_j-x_{k+1} \; ,$$
  and since $h(id,w)>j$, we have $k+1=u(j)$.
\end{proof}

If $\gamma$ is the ascending chain of maximal length \eqref{eq:maxacsendchain} as in Theorem~\ref{thm:main-general}, then
$$u_{k-1}\omega_{i_k} - v\omega_{i_k} = u_{k-1}(\omega_{i_k} - u_{k-1}^{-1}v \omega_{i_k}) = x_{u_{k-1}(i_k)} - x_{v(i_k)} \; ,$$
and therefore
$$
E(\gamma) = \left[ \prod_{(a,b) \in \; \Inv(v)} (x_a-x_b)\right] \prod_{k=1}^m
  \frac{1}{x_{u_{k-1}(i_k)} - x_{v(i_k)}} \; .
$$

We summarize the results of this section and formulate
our version of a type $A$ positive formula for $\tau_u(v)$. Let $C_0(u,v)$ be
the set of ascending chains of maximal length
  \begin{equation}\label{eq:maxacsendchainA}
    \gamma \colon \; u=u_0 \xrightarrow{*(i_1,j_1)} u_1 \xrightarrow{*(i_2,j_2)} u_2 \to \dotsb \to u_{m-1} \xrightarrow{*(i_m,j_m)} u_m = v
  \end{equation}
   that satisfy the condition
  \begin{equation}\label{eq:convexchainsA}
    i_1 \leqslant i_2 \leqslant \dotsb \leqslant i_m \; ,
  \end{equation}
  where $*(i,j)$ means that we are multiplying to the right by the transposition that swaps $i$ and $j$. For such a chain $\gamma$, let
  $$\Inv(v,\gamma) = \Inv(v) \setminus \{ (u_{k-1}(i_k),v(i_k)) \; |
  \; k=1,\ldots,m\}$$
    and
  \begin{equation}\label{eq:Emu0A}
    E(\gamma) =  \prod_{(a,b) \in \; \Inv(v,\gamma)} (x_a-x_b) \; .
  \end{equation}

\begin{thm}\label{thm:main-typeA}{\rm
  The restriction of $\tau_u$ at $v$ is given by
  \begin{equation}\label{eq:posFormulaChains}
    \tau_u(v) = \sum_{\gamma \in C_0(u,v)} E(\gamma) =
    \sum_{\gamma \in C_0(u,v)} \prod_{(a,b) \in \; \Inv(v,\gamma)} (x_a-x_b) \; .
  \end{equation}
}
\end{thm}

\begin{example}
\label{ex:chainContrib} {\rm Consider the permutations
$u=2143$ and $v=3421$ in $S_4$.
There are two chains in $\cC_0(2143,3421)$:
\begin{align*}
  \gamma_1 \; \colon &  \quad 2143 \xrightarrow{*(1,4)} 3142
  \xrightarrow{*(2,3)} 3412 \xrightarrow{*(3,4)} 3421 \\
  \gamma_2 \; \colon &  \quad 2143 \xrightarrow{*(1,4)} 3142 \xrightarrow{*(2,4)}
  3241 \xrightarrow{*(2,3)} 3421 \; .
\end{align*}

For the first chain $u_0=2143$, $u_1=3142$, $u_2=3412$, and
$u_3=v=3421$. The edge $2143 \xrightarrow{*(1,4)} 3142$ deletes the
factor $x_{u_0(1)} - x_{v(1)} = x_2-x_3$. The edge $3142
  \xrightarrow{*(2,3)} 3412$ deletes the factor $x_{u_1(2)} -
  x_{v(2)}=x_1-x_4$. The edge $3412 \xrightarrow{*(3,4)} 3421$
  deletes the factor $x_{u_2(3)}-x_{v(3)}=x_1-x_2$. Therefore
  \begin{align*}
    E(\gamma_1) = &
    \frac{(x_2-x_3)(x_1-x_3)(x_2-x_4)(x_1-x_4)(x_1-x_2)}{(x_2-x_3)(x_1-x_4)(x_1-x_2)}
    = \nonumber \\
    = & (x_1-x_3)(x_2-x_4) = (\alpha_1+\alpha_2)(\alpha_2+\alpha_3) \; .
  \end{align*}
Similarly
  \begin{align*}
    E(\gamma_2) = &
    \frac{(x_2-x_3)(x_1-x_3)(x_2-x_4)(x_1-x_4)(x_1-x_2)}{(x_2-x_3)(x_1-x_4)(x_2-x_4)}
    = \nonumber \\
    = & (x_1-x_3)(x_1-x_2) = (\alpha_1+\alpha_2)\alpha_1 \; .
  \end{align*}
Hence
\begin{equation*}
  \tau_{2143}(3421) = E(\gamma_1)+E(\gamma_2) =
  (\alpha_1+\alpha_2)(\alpha_1+\alpha_2+\alpha_3)\; .
\end{equation*} }
\end{example}

\section{Chains and Subwords}
\label{sec:chains-subwords}

In this section we give a general construction connecting ascending chains and subwords. This construction is valid for all Weyl groups and orderings of simple roots.

Let $u \prec v$ and let $I=[i_1,\ldots,i_k]$ be a reduced word for $v$. Let $\cA(u,v)$ be the set of all ascending chains
  \begin{equation}\label{eq:ascchain}
    \gamma \colon \; u=u_0 \xrightarrow{*s_{\beta_1}} u_1 \xrightarrow{*s_{\beta_2}} u_2 \to \dotsb \to u_{m-1} \xrightarrow{*s_{\beta_m}} u_m = v \; ,
  \end{equation}
from $u$ to $v$, not necessarily of maximal length, and let $\cC(u,v)$ be the subset of $\cA(u,v)$ consisting of the chains that satisfy $h(\beta_1) \leqslant \dotsb \leqslant h(\beta_m)$. Then
$$\cC(u,v) \cap \Sigma(u,v) = \cC_0(u,v) \subset \cC(u,v) \subset \cA(u,v) \; .$$

 Let $\cS(u,I)$ be the set of subwords of $I$ that generate words for $u$ by deleting the zeroes; then $\cR(u,I)$ is a subset of $\cS(u,v)$. We define a function $F_I \colon \cA(u,v) \to \cS(u,I)$ as follows.

Let $\gamma$ be the ascending chain \eqref{eq:ascchain}. We use the edges of $\gamma$, in reverse order, to delete letters from $I$ and get a subword of $I$ that is a word for $u$. Since $u_{m-1} = vs_{\beta_m} \prec v$, the Strong Exchange Property implies that there exists a $j$ (unique, since $I$ is reduced) such that
$$u_{m-1} = s_{i_1} \dotsb s_{i_{j-1}} s_{i_{j+1}} \dotsb s_{i_k} \; ,$$
 hence $I_{m-1} = [i_1,\ldots,i_{j-1},0,i_{j+1},\ldots,i_k] \in \cS(u_{m-1},I)$ is a subword of $I$ that is a word for $u_{m-1}$. Similarly, $u_{m-2}=u_{m-1}s_{\beta_{m-1}}$ and $\ell(u_{m-2}) < \ell(u_{m-1})$, so there exists a subword of $I_{m-1}$, obtained by deleting exactly one letter from $I_{m-1}$, that is a word for $u_{m-2}$. However, since $\ell(u_{m-1})$ is not necessarily equal to $\ell(v)\!-\!1$, the word $I_{m-1}$ may be non-reduced, and in this case the uniqueness of the deleted letter is not guaranteed. We choose the subword for which the deleted letter is the rightmost choice, and get a subword $I_{m-2}$ of $I_{m-1}$ (hence of $I$), that is a word, not necessarily reduced, for $u_{m-2}$. Continuing this process, with the same rule for making a choice, if needed, we get a sequence $I_{m-1}, \ldots,I_0$ of subwords of $I$ such that $I_k \in \cS(u_k,I)$ and $I_k$ is obtained from $I_{k+1}$ by deleting one letter. We define $F_I(\gamma) = I_0$, the last subword in the sequence.

\begin{example}{\rm
Let $u=s_1$ and $v=s_1s_2s_1s_2$ in the Weyl group of the root system $C_2$. Let $I=[1,2,1,2]$ be a reduced word for $v$ and let
$$  \gamma_1 :  \quad s_1 \xrightarrow{*s_{\alpha_1+\alpha_2}} s_2s_1
                         \xrightarrow{*s_{\alpha_1+2\alpha_2}} s_1s_2s_1
                         \xrightarrow{*s_{\alpha_2}} s_1s_2s_1s_2 \; ,$$
be an ascending chain of maximal length from $u$ to $v$. Since $s_{\alpha_1+\alpha_2} = s_1s_2s_1$ and $s_{\alpha_1+2\alpha_2} = s_2s_1s_2$, the subword $F_I(\gamma_1)$ is computed as follows:
\begin{align*}
  [1,2,1,2] \cdot & [2] =  [1,2,1,0] = I_2 \\
  [1,2,1,0] \cdot & [2,1,2] =  [0,2,1,0] = I_1 \\
  [0,2,1,0] \cdot & [1,2,1] =  [0,0,1,0] = I_0
\end{align*}
hence $F_I(\gamma_1) = [0,0,1,0]$. The chains in $\Sigma(u,v)$ are given in \eqref{eq:exampleC2}, and $F_I(\gamma_1) = F_I(\gamma_2)= F_I(\gamma_3) = [0,0,1,0]$ and $F_I(\gamma_4)= [1,0,0,0]$.
}
\end{example}

Note that if $\gamma \in \Sigma(u,v)$, then $\ell(u_{k-1}) = \ell(u_k)-1$ for all $k$, hence the words $I_{m-1}$, \ldots, $I_1$, $I_0$ are reduced. Therefore at each stage there is only one possibility for the deleted letter and no choice is necessary.

\begin{lemma}\label{lem:subw-chains}{\rm
  Let $\gamma \in \cA(u,v)$ be an ascending chain. Then
  $$F_I(\gamma) \in \cR(u,I) \Longleftrightarrow \gamma \in \Sigma(u,v) \; .$$
  Moreover, if $F_I \colon \cC(u,v) \to \cS(u,I)$ is bijective, then $F_I \colon \cC_0(u,v) \to \cR(u,I)$ is also bijective.
  }
\end{lemma}

\begin{proof}
  Let $\gamma \in \cA(u,v)$. Then $F_I(\gamma) \in \cR(u,I)$ if and only if $F_I(\gamma)$ has $\ell(u)$ non-zero letters. This happens if and only if we delete $\ell(v) - \ell(u)$ letters from $I$, and since every deleted letter corresponds to an edge of $\gamma$, if and only if $\gamma$ has $\ell(v) - \ell(u)$ edges, in other words, if and only if $\gamma$ is an ascending chain of maximal length. If $\gamma \in \cC_0(u,v)$ then $\gamma \in \Sigma(u,v)$, hence $F_I(\gamma) \in \cR(u,I)$. Therefore $F_I \colon \cC_0(u,v) \to \cR(u,I)$ is well-defined.

  Suppose that $F_I \colon \cC(u,v) \to \cS(u,I)$ is bijective. Then the restriction of $F_I$ to $\cC_0(u,v)$ is injective. Moreover, if $J \in \cR(u,I)$, then there exists $\gamma \in \cC(u,v)$ such that $F_I(\gamma) = J$. But then $\gamma \in \cC(u,v) \cap \Sigma(u,v) = \cC_0(u,v)$, hence $F_I \colon \cC_0(u,v) \to \cR(u,I)$ is also surjective.
\end{proof}

\section{Equivalence of Positive Formulas}
\label{sec:equivalence}

We show that in the case $G=SL_n(\CC)$, the positive formula \eqref{eq:posFormulaChains} is an alternative formulation of Billey's formula \eqref{eq:posFormula}, by proving that, for a particular choice of a reduced word $I$ for $v$, the function $F_I$ restricts to a bijection $F_I \colon \cC_0(u,v) \to \cR(u,I)$ and $SC(F_I(\gamma),I) = E(\gamma)$ for all $\gamma \in \cC_0(u,v)$.

Let $I(v)$ be the reduced word for $v$ constructed inductively using \eqref{eq:decompA}. Then
$$
I(v)=[I_1,I_2,\ldots, I_{n-1}] \; ,
$$
such that, for all $j=1,..,n-1$,
\begin{itemize}
\item $I_j = [k_j, k_j\!-\!1, \ldots, j\!+\!1,j\,]$ or $I_j = [\,]$, and
\item $h(s_{I_1}\dotsb s_{I_j},v) = h(id, s_{I_j}\dotsb s_{I_1}v) = h(id, s_{I_{j+1}} \dotsb s_{I_{n-1}}) > j$.
\end{itemize}

\begin{example}{\em
  If $v=s_2s_1s_3s_2s_3$, then $I_1 = [2,1]$ , $I_2=[3,2]$
  and $I_3=[3]$. }
\end{example}

We can now prove the equivalence of our positive formula (using chains) with Billey's positive formula (using
subwords).

\begin{thm}\label{thm:equivalence}{\em
  Let $u,v \in S_n$ and let $I=I(v)$. Then  $F_I \colon \cC_0(u,v) \to
  \cR(u,I)$ is a bijection and
  \begin{equation*}
    SC(F_I(\gamma),I) = E(\gamma)
  \end{equation*}
  for all chains $\gamma \in \cC_0(u,v)$.}
\end{thm}

\begin{proof}
  In \cite{Za} we proved that $F_I \colon \cC(u,v) \to \cS(u,I)$ is a bijection. Then, by Lemma~\ref{lem:subw-chains}, the map $F_I \colon \cC_0(u,v) \to \cR(u,I)$ is also a bijection.

  To prove that $SC(F_I(\gamma),I) = E(\gamma)$, we show that the factor of $\Lambda_v^-$ canceled by an edge of $\gamma$ is the same as the factor canceled by the corresponding deleted letter from $I$.

  In \cite{Za} we also proved that if
  \begin{equation}
    \gamma \colon \; u=u_0 \xrightarrow{*s_{\beta_1}} u_1 \xrightarrow{*s_{\beta_2}} u_2 \to \dotsb \to u_{m-1} \xrightarrow{*s_{\beta_m}} u_m = v \; ,
  \end{equation}
is a chain in $\cC_0(u,v)$, then:
\begin{enumerate}
  \item The edge $u_{k-1} \xrightarrow{*s_{\beta_k}} u_k$ deletes a letter from $I_{h(\beta_k)}$, and
  \item The letters deleted in $I_j$ are deleted from left to right.
\end{enumerate}

Let $u_{k-1} \to u_k =u_{k-1}s_{\beta_k}$ be an edge of $\gamma$ and let $j = h(\beta_k)$.

Suppose that $h(\beta_k) = \ldots = h(\beta_{k+q-1}) = j < h(\beta_{k+q})$, hence $\beta_k$ is the $q^{\text{th}}$ occurrence of $j$ from right to left. Let
$$J_{k} = [I_1,\ldots,I_{j-1},I_j',I_{j+1}',\ldots,I_{n-1}'] \in \cR(u_k,I(v))$$
be the subword obtained after using $s_{\beta_m}$, \ldots, $s_{\beta_{k+1}}$ to delete letters from $I(v)$, and let
$$J_{k-1} = [I_1,\ldots,I_{j-1},I_j'',I_{j+1}',\ldots,I_{n-1}'] \in \cR(u_{k-1},I(v))$$
be the subword obtained after using $s_{\beta_{k}}$ to delete one more letter from $J_k$. Then $I_j'$ is of the form
$$I_j' = [k_j, \ldots, \widehat{j_1}, \ldots, \widehat{j_{q-1}}, \ldots, j] \subset I_j \; ,$$
and
$$I_j'' = [k_j, \ldots, \widehat{j_1}, \ldots, \widehat{j_{q-1}}, \ldots, \widehat{j_q}, \ldots, j] \subset I_j' \subset I_j\; ,$$
with $j_1 > j_2 > \dotsb > j_q$, since the letters in $I_j$ are deleted from left to right.

The factor of $\Lambda_v^-$ canceled by the edge $u_{k-1} \to u_k$ is
\begin{align*}
  u_{k-1}\omega_j - v\omega_j = & (u_{k-1}\omega_j - u_k\omega_j) + \dotsb + (u_{k+q-2}\omega_j - u_{k+q-1}\omega_j) = \\
  = & \langle \omega_j, \beta_k \rangle u_{k-1}\beta_k + \dotsb + \langle \omega_j, \beta_{k+q-1} \rangle u_{k+q-2}\beta_{k+q-1} = \\
  = & s_{I_1}\dotsb s_{I_{j-1}} s_{k_j}\dotsb s_{j_1+1}( \widehat{s_{j_1}} \dotsb \widehat{s_{j_2}} \dotsb \widehat{s_{j_q}}\alpha_{j_q} + \dotsb + \alpha_{j_1}) \; .
\end{align*}

But if $j_1 > j_2$, then
$$\widehat{s_{j_1}}\dotsb \widehat{s_{j_2}}\alpha_{j_2}+  \alpha_{j_{\,1}}= s_{j_1}\dotsb \widehat{s_{j_2}}\alpha_{j_2} \; , $$
and then, by induction on $q$,
$$\widehat{s_{j_1}} \dotsb \widehat{s_{j_2}} \dotsb \widehat{s_{j_q}}\alpha_{j_q} + \dotsb + \alpha_{j_1} = s_{j_1}\dotsb s_{j_{q-1}} \dotsb \widehat{s_{j_q}}\alpha_{j_q} \; .$$
Therefore
$$u_{k-1}\omega_j - v\omega_j = s_{I_1}\dotsb s_{I_{j-1}} s_{k_j}\dotsb s_{j_q+1}\alpha_{j_q}\; ,$$
and that is precisely the factor of $\Lambda_v^-$ canceled by the missing letter $j_q$, letter that has been deleted by the edge $u_{k-1} \to u_k$.

This finishes the proof of Theorem~\ref{thm:equivalence}, which was the last piece in the proof of Theorem~\ref{thm:main-general}. We have therefore shown that for $G=SL_n(\CC)$, our positive formula (using chains) is an alternative version of Billey's positive formula (using subwords).
\end{proof}

\end{document}